\documentclass[12pt,letterpaper]{amsart}
\usepackage{amsmath,amssymb,amsfonts,amsthm}
\usepackage{mathrsfs}
\usepackage[all]{xy}
\usepackage{stmaryrd}
\usepackage{mathtools}
\usepackage{indentfirst}
\usepackage{verbatim}
\usepackage{hyperref}
\usepackage{color}
\usepackage[usenames,dvipsnames]{xcolor}
\usepackage{graphicx}
\graphicspath{ {images/} }
\usepackage{tikz}
\usepackage{nicefrac}
\usepackage{xfrac}
\usepackage[shortlabels]{enumitem}

\title[Boundedness of MPDVR of canonical surface foliations]{Boundedness of minimal partial du Val resolutions of canonical surface foliations}
\author{Yen-An Chen}
\thanks{The author was partially supported by NSF research grants no: DMS-1801851, DMS-1840190 and by a grant from the Simons Foundation; Award Number: 256202. }
\address{Department of Mathematics, University of Utah, Salt Lake City, UT 84112, USA}
\email{yachen@math.utah.edu}

\usepackage{csquotes}

\newtheorem*{claim}{Claim}
\newtheorem{thm}{Theorem}[section]
\newtheorem{prop}[thm]{Proposition}

\newtheorem{lem}[thm]{Lemma}

\theoremstyle{definition}
\newtheorem{defn}[thm]{Definition}
\newtheorem*{acks}{Acknowledgements}
\theoremstyle{remark}
\newtheorem{rmk}[thm]{Remark}
\newtheorem*{pf}{Proof}

\newcommand\cA{{\mathcal{A}}}

\newcommand\cE{{\mathcal{E}}}
\newcommand\cF{{\mathcal{F}}}

\newcommand\cL{{\mathcal{L}}}

\newcommand\cO{{\mathcal{O}}}

\newcommand\cQ{{\mathcal{Q}}}

\newcommand\cS{{\mathcal{S}}}

\newcommand\cU{{\mathcal{U}}}

\newcommand\cX{{\mathcal{X}}}
\newcommand\cY{{\mathcal{Y}}}
\newcommand\cZ{{\mathcal{Z}}}

\newcommand\bC{{\mathbb C}}

\newcommand\bN{{\mathbb N}}
\newcommand\bP{{\mathbb P}}
\newcommand\bQ{{\mathbb Q}}
\newcommand\bR{{\mathbb R}}

\newcommand\bZ{{\mathbb Z}}

\newcommand\sF{{\mathscr{F}}}
\newcommand\sG{{\mathscr{G}}}

\newcommand\wt[1]{\widetilde{#1}}

\begin{document}

\maketitle

\begin{abstract}
In this paper, we prove the boundedness of foliated surfaces $(X,\sF)$ which are minimal partial du Val resolutions of canonical models $(X_c,\sF_c)$ of general type. 
For applications, we show the boundedness of non-cusp singularities on canonical models of foliated surfaces of general type and the effective generation on the complement of the cusp singularities. 
\end{abstract}

\section{Introduction}
Recently, there has been significant progress in the study of foliated surfaces of general type. 
By work of McQuillan and Brunella (See \cite{mcquillan2008canonical} and \cite{brunella2015birational}), it is known that smooth foliated surfaces of general type with only canonical singularities admit a unique canonical model. 
It is then natural to wonder if these canonical models have a good moduli theory and, particularly, if they admit a moduli functor. 

A first step in this direction was achieved in \cite{hacon2021birational} where it is shown that foliated surfaces of general type with fixed Hilbert function (of the canonical model) are birationally bounded. 

In this paper, we improve the results of \cite{hacon2021birational} by showing that the canonical models and their minimal partial du Val resolutions are bounded. (See Theorem~\ref{key_thm}.)

\begin{thm}
Fix a function $P : \bZ_{\geq 0} \to \bZ$. 
Then the set $\cS_P$ (see Definition~\ref{S_p}) is bounded. 
More precisely, there exists a projective morphism $\mu : \cX \rightarrow T$ and a family of foliations $\cF$ on $\cX$ over $T$ (see Definition~\ref{family_fol}) where $\cX$ and $T$ are quasi-projective varieties of finite type such that, for any $(X,\sF)\in\cS_P$, there is a $t\in T$ and an isomorphism $\rho : X\rightarrow \cX_t$ such that $\sF \cong \rho^*(\cF\vert_{\cX_t})$.
\end{thm}

For applications, we first improve \cite[Proposition 4.1]{hacon2021birational}. (See Proposition~\ref{dihedral_bdd}.)
\begin{prop}
Fix a function $P:\bZ_{\geq 0} \rightarrow \bZ$. 
There exist constants $C_1$ and $C_2$ such that, for any canonical model $(X_c,\sF_c)$ with Hilbert function $\chi(X_c,mK_{\sF_c})=P(m)$ for $m\geq 0$, the index (resp. embedding dimension) at each dihedral singularity 
is at most $C_1$ (resp. $C_2$). 
\end{prop}

We also answer \cite[Conjecture 1]{hacon2021birational}. (See Theorem~\ref{bir}.)
\begin{thm}
Fix a function $P: \bZ_{\geq 0} \rightarrow \bZ$. 
There exists an integer $m_P$ such that if $(X,\sF)$ is a canonical model of a surface with $\kappa(K_\sF)=2$ and $\chi(mK_\sF)=P(m)$ for all $m\geq 0$, then for any $m>0$ divisible by $m_P$, $|mK_\sF|$ defines a birational map which is an isomorphism on the complement of the cusp singularities. 
\end{thm}

\begin{acks}
The author would like to thank Adrian Langer for his precious feedbacks and Christopher D. Hacon for invaluable discussions and encouragement. 
\end{acks}

\section{Preliminaries}
We will always work over $\bC$. 
For any sheaf $\sF$ on $X$, we denote by $\sF^*$ the dual sheaf $\mathscr{H}\!om_{\cO_X}(\sF,\cO_X)$. 

\subsection{Intersection theory on normal complete surfaces}
We will use Mumford's intersection theory. 
We include the definition for the reader's convenience. 

Let $X$ be a normal complete surface and $f: Y\rightarrow X$ be a proper birational morphism from a smooth surface $Y$, for instance, the minimal resolution of singularities. 
Let $E = \cup_jE_j$ be the exceptional divisor of $f$. 
Since the intersection matrix $(E_i\cdot E_j)$ is negative definite, for any Weil $\bR$-divisor $D$ on $X$, we can define $f^*D$ as $\wt{D} + \sum_ja_jE_j$ where $\wt{D} = f_*^{-1}D$ is the proper transform of $D$ and the $a_j$'s are real numbers uniquely determined by the property that $(\wt{D} + \sum_ja_jE_j)\cdot E_i = 0$ for every $i$. 

For any two Weil $\bR$-divisors $D_1$ and $D_2$ on $X$, we define the intersection number $D_1\cdot D_2 := (f^*D_1)\cdot (f^*D_2)$. 
Note that if $D_1$ is Cartier and $D_2$ is an irreducible curve, then $D_1\cdot D_2$ agrees with the degree of $\cO_{D_2}(D_1)$. 

\subsection{Foliations on normal varieties}
\begin{defn}
Let $X$ be a normal variety over $\bC$. 
A \emph{foliation} $\sF$ of rank one on $X$ is a torsion-free rank one quotient sheaf $\Omega_{\sF}$ of $\Omega_X^{**}$. 
\end{defn}

\begin{lem}\label{normal_fol_equiv}
Let $X$ be a normal variety over $\bC$. 
A foliation $\sF$ of rank one $X$ is equivalent to a saturated rank one coherent subsheaf $T_\sF$ of the tangent sheaf $T_X$. 
Here a subsheaf $T_\sF$ of $T_X$ is called saturated if $T_X/T_\sF$ is torsion-free. 
\end{lem}
\begin{pf}
Given a foliation $\sF$ on $X$, we have an exact sequence 
\[\xymatrix{0\ar[r] & K \ar[r] & \Omega_X^{**} \ar[r]^-{\rho_\sF} & \Omega_\sF \ar[r] & 0}\]
where $K$ is the kernel of $\rho_\sF$. 
Dualizing this sequence, we get 
\[\xymatrix{0 \ar[r] & \Omega_\sF^* \ar[r]^-{\rho_\sF^*} & T_X \ar[r] & A \ar[r] & 0}\]
where $A$ is the cokernel of $\rho_\sF^*$. 
Note that $A$ is a subsheaf of $K^*$ and $K^*$ is torsion-free. 
So $A$ is also torsion-free, and thus $\Omega_\sF^*$ is saturated. 
Therefore, this gives a function $F$ from the set of foliations of rank one on $X$ to the set of saturated rank one coherent subsheaves of the tangent sheaf $T_X$. 

Conversely, given a saturated coherent subsheaf $T_\sF$ of $T_X$, we have an exact sequence 
\[\xymatrix{0 \ar[r] & T_\sF \ar[r] & T_X \ar[r]^-\pi & T_X/T_\sF \ar[r] & 0.}\] 
Dualizing the sequence, we have 
\[\xymatrix{0\ar[r] & (T_X/T_\sF)^* \ar[r]^-{\pi^*} & \Omega_X^{**} \ar[r] & B \ar[r] & 0}\]
where $B$ is the cokernel of $\pi^*$. 
Notice that $B$ is a subsheaf of $T_\sF^*$ and $T_\sF^*$ is torsion-free. 
Thus $B$ is also torsion-free and therefore gives a foliation on $X$. 
Hence, this gives a function $G$ from the set of saturated rank one coherent subsheaves of the tangent sheaf $T_X$ to the set of foliations of rank one on $X$. 

Notice that both $G\circ F$ and $F\circ G$ are generically identity functions. 
\begin{claim}
They are indeed identity functions up to isomorphism. 
\end{claim}
\begin{pf}
Notice that we have the following exact sequences: 
\[\xymatrix{0 \ar[r] & K \ar[r] & \Omega_X^{**} \ar[r] \ar@{=}[d] & \Omega_\sF \ar[r] & 0 \\
0 \ar[r] & A^* \ar[r] & \Omega_X^{**} \ar[r] & G(F(\Omega_\sF)) \ar[r] & 0.}\]
Since $G(F(\Omega_\sF))$ is torsion-free and the map 
\[\xymatrix{K\ar[r] & \Omega_X^{**} \ar[r] & G(F(\Omega_\sF))}\] 
is generically zero, the image is in fact zero in $G(F(\Omega_\sF))$. 
Thus, we have $K\subset A^*$. 
Similarly, we have $A^* \subset K$. 
Therefore, we have $K = A^*$ and hence $\Omega_\sF \cong G(F(\Omega_\sF))$. 

Similarly, we could show that $T_\sF = F(G(T_\sF))$. 
\end{pf}
\qed
\end{pf}

\begin{defn}\label{family_fol}
Given a family $\mu:\cX\rightarrow T$ where $\cX$ is normal. 
A foliation $\cF$ (of rank one) on $\cX$ is called a \emph{family of foliations} over $T$ if $\Omega_{\cX}^{**} \rightarrow \Omega_\cF$ factors through $\Omega_{\cX/T}^{**}$.
\end{defn}

\subsection{Foliations on normal surfaces}
Now we focus on the case when $X$ is a normal surface. 
By Lemma~\ref{normal_fol_equiv}, a foliation $\sF$ on a normal surface $X$ is given by a saturated subsheaf $T_\sF$ of the tangent sheaf $T_X$. 

A  point $p$ on $X$ is called a \emph{singular} point of the foliation $\sF$ if either it is a singular point of $X$ or a point at which the quotient $T_X/T_\sF$ is \emph{not} locally free. 
Since $X$ is normal, we have only isolated singularities for $\sF$. 
\begin{defn}\label{defn_pair}
A \emph{foliated surface} is a pair $(X,\sF)$ consisting of a surface $X$ and a foliation $\sF$ on $X$. 
\end{defn}

Since $T_X = \Omega_X^*$ is reflexive, we have that $T_\sF$ is also reflexive because it is a saturated subsheaf of a reflexive sheaf $T_X$. 
Hence, we define the \emph{canonical divisor} $K_\sF$ of the foliation as a Weil divisor on $X$ with $\cO_X(-K_\sF)\cong T_\sF$. 

\begin{defn}
Let $(X,\sF)$ be a foliated surface. 
Given any proper biration morphism of normal surfaces $f : Y\rightarrow X$ and a foliation $\sF$ on $X$, then we define the \emph{pullback foliation} $f^*\sF$ as follows: 

Let $U$ be the largest open subset such that $V:= f^{-1}(U) \rightarrow U$ is an isomorphism. 
Note that $\sF\vert_U \subset T_U \cong T_V$. 
By \cite[Exercise II.5.15]{hartshorne1977algebraic}, we have a coherent subsheaf $\sG$ of $T_Y$ such that $\sG\vert_V = \sF\vert_U \subset T_V $. 
Then the pullback foliation $f^*\sF$ is defined to be the saturation of $\sG$. 
By \cite[Lemma 1.8]{hacon2021birational}, this definition is well-defined. 

Also if $\sG$ is a foliation on $Y$, then we can define the \emph{pushforward foliation} $f_*\sG$ by taking the saturation of the image of the composition 
\[\xymatrix{f_*T_\sG\ar[r] & f_*T_Y\ar[r] & (f_*T_Y)^{**} = T_X.}\]
\end{defn}

\begin{defn}
Let $(X,\sF)$ be a foliated surface and let $f : Y\rightarrow X$ be a proper birational morphism. 
For any divisor $E$ on $Y$, we define the \emph{discrepancy} of $\sF$ along $E$ to be $a_E(\sF) = \textnormal{ord}_E(K_{f^*\sF}-f^*K_\sF)$. 
We say  
\[x \textnormal{ is } \left\{\begin{array}{ll}
\textnormal{terminal} & \textnormal{ if } a_E(\sF)> 0 \\
\textnormal{canonical} & \textnormal{ if } a_E(\sF)\geq 0 
\end{array}\right. 
\textnormal{ for every divisor $E$ over $x$.} \]
\end{defn}

\begin{defn}
We define the \emph{Kodaira dimension} of $\sF$ as 
\[\kappa(\sF) = \max\{\dim\phi_{mK_\sF}(X)\vert\,m\in\bN_{>0} \textnormal{ and } h^0(X,mK_\sF)>0\}\] 
where by convention $\kappa(\sF)=-\infty$ if the set above is empty, that is, $h^0(X,mK_\sF)=0$ for all $m\in\bN_{>0}$. 
We say that $\sF$ is \emph{of general type} if $\kappa(K_\sF) = \dim X$. 
\end{defn}

\begin{defn}\label{index}
Let $(X,\sF)$ be a foliated surface. 
We have the following definitions for indices.
\begin{enumerate}
\item The index $i(X)$ of $X$ is the smallest positive integer $m$ such that $mD$ is Cartier for every Weil divisor $D$ on $X$. 
(We set $i(X)=\infty$ if there is no such $m$.)
\item The index $i(K_X)$ of $K_X$ is the smallest positive integer $m$ such that $mK_X$ is Cartier. 
(We set $i(K_X)=\infty$ if $K_X$ is not $\bQ$-Cartier.)
\item The index $i(\sF)$ of a foliation $\sF$ on $X$ is the smallest positive integer $m$ such that $mK_\sF$ is Cartier. 
(We set $i(\sF)=\infty$ if $K_\sF$ is not $\bQ$-Cartier.)
\item The $\bQ$-index $i_\bQ(\sF)$ of a foliation $\sF$ on $X$ is the smallest positive integer $m$ such that $mK_\sF$ is Cartier for all $\bQ$-Gorenstein points of the foliation. 
\end{enumerate}
\end{defn}

\begin{defn}[Canonical model]
A foliated surface $(X,\sF)$ is called a \emph{canonical model} if $\sF$ is a foliation with only canonical singularities on a normal surface $X$, $K_\sF$ is nef, and $K_\sF\cdot C = 0$ implies $C^2\geq 0$ for all irreducible curves $C$. 
\end{defn}

\begin{lem}[{\cite[Lemma 1.11]{hacon2021birational}}]\label{ample}
Let $(X,\sF)$ be a canonical model. 
If the foliation $\sF$ is of general type, then $K_{\sF}$ is numerically ample, that is $K_{\sF}^2>0$ and $K_{\sF}\cdot C>0$ for any irreducible curve $C$ on $X$. 
\end{lem} 

In \cite{hacon2021birational}, Hacon and Langer show that for a canonical model $(X,\sF)$ with the fixed Hilbert function $\chi(X,mK_\sF)$, many invariants are bounded. 
Precisely, 
\begin{prop}[{\cite[Proposition 4.1]{hacon2021birational}}]\label{bdd}
Fix a function $P:\bZ_{\geq 0} \rightarrow \bZ$. 
Then there exist some constants $B_1$, $B_2$, $B_3$, and $B_4$ such that, for any canonical model $(X,\sF)$ with $\chi(X,mK_\sF) = P(m)$, the intersection numbers $K_\sF^2 = B_1$, $K_\sF\cdot K_X = B_2$, $\chi(\cO_X) = B_3$, and the number of cusps of $X$ is equal to $B_4$. 
Moreover, there exists some constants $C_1$ and $C_2$ such that the number of terminal and dihedral singularities of $(X,\sF)$ is at most $C_1$ and the index of $X$ at any terminal foliation singularity is at most $C_2$. 
In particular, we have $i_\bQ(\sF)\leq 2C_2$.
\end{prop}

\subsection{A criterion for very-ampleness}
\begin{thm}[{\cite[Theorem 0.2]{langer2001adjoint}}]\label{va_eff}
Let $X$ be a normal projective surface and $M$ be a $\bQ$-divisor on $X$ such that $K_X+\lceil M\rceil$ is Cartier. 
Let $\zeta$ be a 0-dimensional subscheme of $X$ and $\delta_\zeta>0$ be a certain number associated to $\zeta$. 
If $M^2>\delta_\zeta$ and $M\cdot C\geq \frac{1}{2}\delta_\zeta$ for every curve $C$ on $X$, then 
\[\xymatrix{H^0(X,\cO_X(K_X+\lceil M\rceil)) \ar[r] & \cO_\zeta(K_X+\lceil M\rceil)}\] 
is surjective. 
\end{thm}

\begin{rmk}\label{delta}
We consider the case when $|\zeta| = 2$. 
\begin{enumerate}
\item If the support of $\zeta = \{x_1$, $x_2\}$ is two distinct points, then $\delta_{\zeta} = \delta_{x_1}+\delta_{x_2}$. (\cite[1.1.3]{langer2001adjoint})
\item $\delta_x \leq 4$ if $x$ is a smooth point or a du Val singularity. (\cite[0.3.2]{langer2001adjoint})
\item If $\zeta$ is supported at only one point $x$, then 
\begin{enumerate}
\item $\delta_\zeta \leq 8$ if $x$ is a smooth point or a du Val singularity. (\cite[0.3.2]{langer2001adjoint})
\item $\delta_\zeta \leq 4(\textnormal{edim}_xX-1)$ if $x$ is a rational (non-smooth) singularity where $\textnormal{edim}_xX$ is the embedding dimension of $X$ at $x$. (\cite[0.3.4]{langer2001adjoint})
\end{enumerate}
\end{enumerate}
\end{rmk}

\begin{prop}[{\cite[Proposition II.7.3]{hartshorne1977algebraic}}]\label{sep}
Let $X$ be a projective scheme over $k$ and $\cL$ be a line bundle on $X$. 
Then $\cL$ is very ample if and only if it separates points and tangents.
\end{prop}

\subsection{Riemann-Roch theorem on normal surfaces}
\begin{thm}[{\cite{reid1987canonical},\cite[Section 3]{langer2000chern}}]\label{RR}
Let $X$ be a normal projective surface and $D$ be a Weil divisor on $X$. 
Then we have 
\[\chi(D) = \frac{1}{2}D\cdot (D-K_X) + \chi(\cO_X)+\sum_{x\in\textnormal{Sing}(X)} a(x,D)\]
where $a(x,D)$ is a local contribution of $\cO_X(D)$ at $x$ depending only on the local isomorphism class of the reflexive sheaf $\cO_X(D)$ at $x$.
\end{thm}

If $X$ has only quotient singularities, then the theorem follows from \cite{reid1987canonical}. 
In general, the theorem follows from \cite[Section 3]{langer2000chern}.

\subsection{Quot scheme}
In this section, we recall some properties of Quot schemes. 

\begin{defn}
Let $X \rightarrow S$ be a scheme of finite type over a noetherian base scheme $S$ and $\cE$ be a coherent sheaf on $X$. 
We consider the functor 
\[\xymatrix{\cQ uot_{\cE/X/S} : (\textnormal{Sch}/S)^{\textnormal{op}} \ar[r] & (\textnormal{Sets})}\]
sending $T\rightarrow S$ to the set 
\[\cQ uot_{\cE/X/S}(T) = \left\{(\cF,q)\Big\vert 
\begin{array}{l}
\cF\in\textnormal{Coh}(X_T), \textnormal{Supp}(\cF) \textnormal{ is proper over } T \\
\cF \textnormal{ is flat over } T, q:\cE_T\rightarrow \cF \textnormal{ is surjective}
\end{array}\right\}\Big/\sim
\]
where $X_T = X\times_ST$, $\cE_T = \textnormal{pr}_X^*\cE$ under the projection $\textnormal{pr}_X : X_T \rightarrow X$, and $(\cF,q) \sim (\cF', q')$ if $\ker q = \ker q'$. 
\end{defn}

\begin{defn}\label{Hilbert_poly}
Let $X \rightarrow S$ be a scheme of a finite type over a noetherian base scheme $S$ and $\cA$ be a relatively ample line bundle. 
For any coherent sheaf $\cF$ on $X$ flat over $S$, we define the Hilbert polynomial $\Phi_\cF$ with respect to $\cA$ as \[\Phi_\cF(m) = \chi(\cF_s(m)) = \sum_{i=0}^{\dim X}(-1)^ih^i(X_s, \cF_s\otimes\cA_s^{\otimes m})\]
for any fixed $s\in S$. 
Note that this definition is independent of the choice of $s\in S$. 
\end{defn}

\begin{prop}
Fix any relatively ample line bundle $\cA$ on $X$. 
There is a natural stratification: 
\[\cQ uot_{\cE/X/S} = \coprod_{\Phi\in\bQ[\lambda]}\cQ uot_{\cE/X/S}^{\Phi,\,\cA}\]
where $\cQ uot_{\cE/X/S}^{\Phi,\,\cA}(T) = \{(\cF,q)\in\cQ uot_{\cE/X/S}(T)\,\vert\, \Phi_\cF = \Phi\}$. 
\end{prop}

\begin{thm}[Grothendieck]\label{quot_exist}
The functor $\cQ uot_{\cE/X/S}^{\Phi,\,\cA}$ is representable by a projective scheme $\textnormal{Quot}_{\cE/X/S}^{\Phi,\,\cA}$ over $S$.
\end{thm}

\section{Minimal partial du Val resolution}
Fix an integral-valued function $P : \bZ_{\geq 0} \rightarrow \bZ$. 
Given any canonical model $(X_c,\sF_c)$ of general type with Hilbert function $\chi(mK_{\sF_c}) = P(m)$ for all $m\in\bZ_{\geq 0}$. 
Let $X^m$ be the minimal resolution of $X_c$ at points at which $\sF_c$ is canonical \emph{non-terminal}, and $\sF^m$ be the pullback foliation. 
We denote the associated morphism as 
\[\xymatrix{g : (X^m,\sF^m) \ar[r] & (X_c,\sF_c).}\] 
By \cite[Theorem III.3.2]{mcquillan2008canonical}, the connected components of the exceptional divisors $E_j$'s
belong to one of the following types: 
\begin{enumerate}
\item A chain of smooth rational curves whose dual graph is of $A_n$ type. 
More precisely, it consists of either two $(-1)$-$\sF^m$-curves of self-intersection $-2$ joined by a bad tail or a chain of $(-2)$-$\sF^m$-curves.
\item Two $(-1)$-$\sF^m$-curves of self-intersection $-2$ joined by a bad tail which itself connects to a chain of $(-2)$-$\sF^m$-curves. Its dual graph is of $D_n$ type. 
\item Elliptic Gorenstein leaves. That is, either a cycle of $(-2)$-$\sF^m$-curves or a rational curve with only one node.
\end{enumerate}

Let $h: X^m \rightarrow X$ be the (classical) relative canonical model over $X_c$, which is given by a sequence of contractions of smooth rational curves $C$ in the fiber of $g$ with $C^2=-2$, and let $\sF$ be the pushforward foliation via $h$. 
So we have a morphism of foliated surfaces $f: (X,\sF) \rightarrow (X_c,\sF_c)$ such that $g = f\circ h$. 
We summarize all morphisms of foliated surfaces in the following diagram.
\[\xymatrix{(X^m,\sF^m) \ar_-g[dd] \ar^-h[rd] & \\ & (X,\sF) \ar^-f[ld] \\ (X_c,\sF_c) & }\]

We give the following definition. 
\begin{defn}\label{mpcr}
We call $f : (X,\sF) \rightarrow (X_c,\sF_c)$ as above the \emph{minimal partial \mbox{du Val} resolution} of the given canonical foliation $(X_c,\sF_c)$.
\end{defn}

\begin{rmk}
Note that any minimal partial du Val resolution $(X,\sF)$ determines its canonical model $(X_c,\sF_c)$ uniquely. 
\end{rmk}

\begin{defn}\label{S_p}
Fix a function $P : \bZ_{\geq 0} \rightarrow \bZ$. 
Let $\cS_P$ be the set of foliated surfaces $(X,\sF)$ which are minimal partial du Val resolutions of canonical models $(X_c,\sF_c)$ of general type with Hilbert function $\chi(mK_{\sF_c}) = P(m)$. 
\end{defn}

Now, we would like to prove the boundedness of $\cS_P$. 
\begin{thm}\label{key_thm}
Fix a function $P : \bZ_{\geq 0} \to \bZ$. 
Then the set $\cS_P$ is bounded. 
More precisely, there exists a projective morphism $\mu : \cX \rightarrow T$ and a family of foliations $\cF$ on $\cX$ over $T$ (see Definition~\ref{family_fol}) where $\cX$ and $T$ are quasi-projective varieties of finite type such that, for any $(X,\sF)\in\cS_P$, there is a $t\in T$ and an isomorphism $\rho : X\rightarrow \cX_t$ such that $\sF \cong \rho^*(\cF\vert_{\cX_t})$. 
\end{thm}

\begin{pf}
Given any foliated surface $(X,\sF)\in\cS_P$. 
Let $(X_c,\sF_c)$ be the canonical model of $(X,\sF)$ and $f:(X,\sF)\rightarrow (X_c,\sF_c)$ be the associated morphism. 
Note that $f^*K_{\sF_c} = K_\sF$. 
Let $E=\bigcup_iE_i$ be the exceptional divisor of $f$. 
We divide the proof into several steps. 
\begin{enumerate}
\item\label{adjoint_ample} 
\emph{Claim.} $K_X+mi(\sF)K_\sF$ is nef when $m \geq 3$ and is ample when $m\geq 4$. 
\begin{pf}
Given any irreducible curve $C$ on $X$. 
If $C=E_i$ for some $i$, then 
\[(K_X+mi(\sF)K_\sF)\cdot C = K_X\cdot C > 0\]
since $K_X$ is ample over $X_c$. 
If $C$ is not contained in the exceptional divisor, then we have two cases: 
\begin{enumerate}
\item If $K_X\cdot C\geq 0$, then 
\[(K_X+mi(\sF)K_\sF)\cdot C \geq mi(\sF)K_\sF\cdot C = mi(\sF)K_{\sF_c}\cdot f_*C > 0\] 
since $K_{\sF_c}$ is numerically ample by Lemma~\ref{ample}. 

\item If $K_X\cdot C<0$, by \cite[Proposition 3.8]{fujino2012minimal}, every $K_X$-negative extremal ray is spanned by a rational curve $C$ with $0<-K_X\cdot C\leq 3$. 
Thus, if we write $C \equiv \sum a_iC_i+B$ where $B\cdot K_X\geq 0$, the $\bR_{\geq 0}[C_i]$'s are $K_X$-negative extremal rays, the $a_i$'s are non-negative, and at least one of $a_i$ is positive. 
We have seen above that $(K_X+mi(\sF)K_\sF)\cdot B\geq 0$. 
Thus, we have 
\begin{eqnarray*}
(K_X+mi(\sF)K_\sF)\cdot C &\geq & \sum_ia_i \Big(K_X\cdot C_i+mi(\sF)K_\sF\cdot C_i\Big) \\
&\geq & \sum_ia_i(- 3+m) 
\end{eqnarray*}
is non-negative if $m\geq 3$ and positive if $m\geq 4$. 
\end{enumerate} 
When $m\geq 4$, we notice that 
\[K_X+mi(\sF)K_\sF = K_X+3i(\sF)K_\sF + (m-3)i(\sF)K_\sF\] 
is big since it is the sum of a nef divisor and a big divisor. 
\end{pf}

\item \emph{Claim.} $i(K_X)$ and $i(K_\sF)$ are bounded.
\begin{pf}
By Proposition~\ref{bdd}, $i(K_X)$ and $i(K_\sF)$ at points at which $\sF$ is terminal are bounded in terms of $\chi(mK_{\sF_c})$. 
Other foliation singularities are at du Val singularities at which $\sF$ is canonical. 
Then we have $i(K_X)=1$ and $i(K_\sF)\leq 2$ by \cite[Proposition 2.4]{hacon2021birational}.  
\end{pf}

\item \emph{Claim.} The embedding dimension at each point at which $\sF$ is singular is bounded.
\begin{pf}
By Proposition~\ref{bdd}, the index at each point at which $\sF$ is terminal is bounded in terms of $\chi(mK_{\sF_c})$. 
Indeed, the proof of Proposition~\ref{bdd} in~\cite{hacon2021birational} shows that the order of each cyclic group acting at each point at which $\sF$ is terminal is bounded. 
Thus, the embedding dimension at each point at which $\sF$ is terminal is bounded. 
See for example \cite[Corollary 2.5]{reid2012surface}. 
In fact, if the action is $\frac{1}{n}(1,q)$ and $\frac{n}{n-q} = [a_1$, $\ldots$, $a_k]$, then the embedding dimension is at most $k+2$.

Other foliation singularities are at du Val singularities, at which the embedding dimensions are $3$.  

Thus, we have a uniform bound $\delta$ such that $\delta_\zeta\leq\delta$ for $|\zeta|=2$. 
See Remark~\ref{delta}.
\end{pf}

\item We consider the $\bQ$-Cartier divisor
\begin{eqnarray*}
L &=& \alpha i(K_X)(K_X+ 4i(\sF)K_\sF) -K_X \\
&=& (\alpha i(K_X)-1)K_X + 4\alpha i(K_X)i(\sF)K_\sF.
\end{eqnarray*}
where $\alpha$ will be determined later.
Notice that $L$ can be written in the following way:
\[L = (\alpha i(K_X)-1)(K_X+3i(\sF)K_\sF)+(\alpha i(K_X)+3)i(\sF)K_\sF.\] 
Since $i(\sF)K_\sF$ is Cartier, we can fix an $\alpha$, which only depends on the Hilbert function $\chi(mK_{\sF_c})$, such that for any irreducible curve $C$, we have 
\begin{eqnarray*}
L^2 &\geq & (\alpha i(K_X)+3)^2 > \delta \textnormal{ and}\\
L\cdot C &\geq & \alpha i(K_X)+3 > \frac{\delta}{2}. 
\end{eqnarray*}
Thus, by Theorem~\ref{va_eff}, we have 
\[K_X+L=K_X+\lceil L\rceil = \alpha i(K_X)(K_X+ 4i(\sF)K_\sF)\] 
separates points and tangents. 
Hence, by Theorem~\ref{sep}, $K_X+L$ is very ample.

\item \emph{Claim.} The set of surfaces $X$ where $(X,\sF)\in\cS_P$ forms a bounded family. 
\begin{pf}
Let $\varphi: X \rightarrow \bP^N$ be an embedding defined via $|K_X+L|$. 
We may assume that $N$ is minimal in the sense that $\varphi(X)$ is not contained in any hyperplane of $\bP^N$. 
It is known (for a reference, see \cite[Proposition 0]{eisenbud1987varieties}) that 
\[\deg\varphi(X)\geq 1+(N-\dim\varphi(X)) = N-1.\] 

By \cite[Lemma 3.5]{hacon2021birational}, we know that $\gamma K_{\sF_c}-K_{X_c}$ is pseudo-effective for 
\[\gamma = \max\left\{\frac{2K_{\sF_c}\cdot K_{X_c}}{K_{\sF_c}^2}+3i_\bQ(\sF_c),0\right\}\] 
which is bounded in terms of $\chi(mK_{\sF_c})$. 
Thus, we have 
\begin{eqnarray*}
\deg\varphi(X) 
&=& \textnormal{vol}(K_X+L) \\
&=& \textnormal{vol}\Big(\alpha i(K_X)(K_X+ 4i(\sF)K_\sF)\Big) \\
&\leq & \textnormal{vol}\Big(\alpha i(K_X)(K_{X_c}+ 4i(\sF)K_{\sF_c})\Big) \\
&\leq & \textnormal{vol}\Big(\alpha i(K_X)(K_{X_c}+ 4i(\sF)K_{\sF_c}) + \alpha i(K_X)(\gamma K_{\sF_c}-K_{X_c}) \Big) \\
&=& \textnormal{vol}\Big(\alpha i(K_X)(\gamma + 4i(\sF))K_{\sF_c} \Big) \\
&=& \Big(\alpha i(K_X)(\gamma + 4i(\sF))K_{\sF_c} \Big)^2
\end{eqnarray*}
which is bounded. 
Therefore, $N$ is also bounded. 

Hence, by the boundedness result about Chow varieties, we have the boundedness of surfaces $X$ where $(X,\sF)\in\cS_P$.
\end{pf}

\item Now, we would like to show the boundedness of foliations. 

We have shown that there exists a projective morphism $\mu : \cX\rightarrow T$ where $\cX$ and $T$ are quasi-projective varieties of finite type  
such that for any $(X,\sF)\in\cS_P$, 
there is a $t\in T$ such that $\cX_t\cong X$ 
and by construction we have 
\[\xymatrix{\cX \ar@{^{(}->}[r] \ar[d]_-\mu & \bP^N_T \ar[ld] \\ T & }\]
and $\cO_X(mK_\sF+nK_X) = \cO_{\cX_t}(1)$ 
where $m=4\alpha i(K_X)i(\sF)$ and $n = \alpha i(K_X)$. 
Note that, by Proposition~\ref{bdd}, $\alpha$, $i(K_X)$, and $i(\sF)$ are fixed and determined by the function $P$. 
Without loss of generality, we may assume that $T$ is the closure of $T'$ where
\[T' = \{t\,\vert\, \cX_t\cong X \textnormal{ for some } (X,\sF)\in \cS_P\}.\] 

\item We may first assume that $T$ is reduced. 
\begin{claim}
We may assume that $T$ is irreducible. 
\end{claim}
\begin{pf}
Let $T = \bigcup_iT_i$ be the decomposition of $T$ into irreducible components. 
If there is a family of foliations $\cF_i$ on each $\cX_i = \cX\times_T T_i$ such that 
for any $(X,\sF)\in\cS_P$, there exist an $i$ and $t\in T_i$ 
with $(X,\sF)\cong (\cX_i,\cF_i)_t$, 
then the disjoint union of these families gives the required family. 
\end{pf}

\begin{claim}
We may assume that $T$ is smooth. 
\end{claim}
\begin{pf}
This is achieved by the existence of a finite stratification of $T$ by smooth locally closed subsets. 
\end{pf}

\item Recall that, by Serre's criterion, $\cX$ is normal if and only if it is $R_1$ and $S_2$. 

\begin{claim}
We may assume that $\cX$ is normal. 
\end{claim}
\begin{pf}
Suppose that $\cX_t$ is $R_1$ (resp. $S_2$), then there exists an open subset $U\subset T$ containing $t$ such that $\cX_U := \cX\times_T U$ is $R_1$ (resp. $S_2$) and every fiber $\cX_t$ is also $R_1$ (resp. $S_2$) for $t\in U$. 
(For a reference, see \cite[Th\'{e}or\`{e}me 12.2.4]{grothendieck1966elements}.) 
It follows that there is an open subset $V \subset T$ such that $\cX_V := \cX\times_T V$ is normal, every fiber $\cX_t$ is normal for $t\in V$, and for all $X\in \cS$, there exists a $t\in V$ such that $X\cong \cX_t$. 
Hence, by shrinking $T$, we may assume that $\cX$ is normal. 
\end{pf}

\item\label{Q_Cartier} \emph{Claim.} We may assume that $K_\cX$ is $\bQ$-Cartier. 
\begin{pf}
This is well-known; see for example \cite[Lemma 1.27]{simpson1994moduli}. 
We also include one proof. 

Since $\cX$ is normal, the canonical sheaf $\cO_\cX(K_\cX)$ is defined as a divisorial sheaf. 
We fix a $t\in T'$ and put $X=\cX_t$. 
Then $i(K_X)K_{\cX_t}$ is Cartier. 
We fix a $d$ such that $\cO_{\cX_t}(i(K_X)K_{\cX_t})(d)$ is generated by global sections. 
Then by shrinking $T$, we may assume that 
\[\mu_*\cO_\cX(i(K_X)K_\cX)(d)\otimes k(t) \rightarrow H^0(\cX_t,\cO_{\cX_t}(i(K_X)K_{\cX_t})(d))\]
is surjective by \cite[Theorem III.12.8 and Corollary III.12.9]{hartshorne1977algebraic}. 
Note that the sheaf $\cO_{\cX_t}(i(K_X)K_{\cX_t})(d)$ is locally generated by some section of 
\[H^0(\cX_t,\cO_{\cX_t}(i(K_X)K_{\cX_t})(d)).\] 
By lifting this section to a section of $H^0(\cX,\cO_\cX(i(K_\cX)K_\cX)(d))$ using the surjection above and by Nakayama's lemma, the sheaf $\cO_\cX(i(K_X)K_\cX)(d)$ is locally generated by one section on a neighborhood of $\cX_t$. 
Therefore, $i(K_X)K_\cX$ is Cartier on a neighborhood of $\cX_t$. 
After shrinking $T$ further, we may assume that $i(K_X)K_\cX$ is Cartier. 
\end{pf}

\item Now we consider the function 
\begin{eqnarray*}
\Phi(\ell) &=& \chi(X,\cO_X(K_\sF+\ell(mK_\sF+nK_X))) \\
&=& \chi(X,\cO_X((m\ell+1)K_\sF+n\ell K_X)) \\
&=& \frac{1}{2}\Big((m\ell+1)K_\sF+n\ell K_X\Big)\cdot\Big((m\ell+1)K_\sF+(n\ell-1)K_X\Big) \\
& & \quad + \chi(\cO_X) + \sum_{x\in\textnormal{Sing}(X)}a\big(x,(m\ell+1)K_\sF+n\ell K_X\big). 
\end{eqnarray*}
Note that $mK_\sF$ and $nK_X$ are Cartier, we have that 
\[a\big(x,(m\ell+1)K_\sF+n\ell K_X\big) = a(x,K_\sF)\] 
for any singularity $x\in X$. 
Also we have
\begin{eqnarray*}
\chi(\cO_X) + \sum_{x\in\textnormal{Sing}(X)}a(x,K_\sF) 
&=& \chi(K_\sF) - \frac{1}{2}K_\sF\cdot(K_\sF-K_X) \\
&=& P(1) - \frac{1}{2}K_\sF\cdot(K_\sF-K_X).
\end{eqnarray*}
Thus, we have that 
\begin{eqnarray*}
\Phi(\ell) &=& \frac{1}{2}\Big((m\ell+1)K_\sF+n\ell K_X\Big)\cdot\Big((m\ell+1)K_\sF+(n\ell-1)K_X\Big) \\
& & \quad + P(1) - \frac{1}{2}K_\sF\cdot(K_\sF-K_X)
\end{eqnarray*}
is a polynomial with coefficients in terms of $P(1)$, $K_\sF^2$, $(K_\sF\cdot K_X)$, and $K_X^2$. 
Since the function $P$ is fixed, we have $K_\sF^2$ and $K_\sF\cdot K_X$ are fixed by Proposition~\ref{bdd}. 
Moreover, $K_X^2$ has only finitely many possibilities since the surfaces $X$ are in a bounded family. 
Therefore, we have finitely many possible polynomials $\Phi(\ell)$. 

Notice that $\frac{m+1}{n-1}>\frac{m}{n} = 4i(\sF)\geq 4$. 
So we have $(m+1)K_\sF+(n-1)K_X$ is ample by~(\ref{adjoint_ample}) and thus $\Phi(1)=h^0(X,\cO_X((m+1)K_\sF)+nK_X)$. 
Let $\varphi$ be the maximum of $\Phi(1)$ among all possible polynomials $\Phi(\ell)$.

\item Applying \cite[Theorem 21]{kollar2008hulls} to the projective morphism $\mu : \cX \rightarrow T$ and the coherent sheaf $\Omega_{\cX/T}$, after shrinking $T$ further, we may assume that $\Omega_{\cX/T}$ has a relative hull $\Omega_{\cX/T}^{[**]} = \Omega_{\cX/T}^{**}$. 
Note that the relative hull behaves well under base change. 

For any polynomial $\Psi(\ell) := \Phi(\ell)-s$ where $0\leq s\leq \varphi$, we consider the Quot scheme 
\[Q := \textnormal{Quot}_{\Omega_{\cX/T}^{**}/\cX/T}^{\Psi(\ell),\,\cO_\cX(1)},\] 
which is projective over $T$. 
Let $\cU := \cX\times_T Q$. 
By a flattening stratification on $Q$, we may assume that the projection morphism $\cU \rightarrow Q$ is flat. 
Also, by the universal property of Quot schemes, there is a (universal) quotient sheaf $\cL$ of $\Omega_{\cU/Q}$ on $\cU$. 
Moreover, by \cite[Th\'{e}or\`{e}me 12.2.1 (v)]{grothendieck1966elements}, the loci 
\[\{q\in Q\, : \,\cL\vert_{\cU_q} \textnormal{ is torsion-free}\} \textnormal{ and } \{q\in Q\, : \,\cL^*\vert_{\cU_q} \textnormal{ is reflexive}\}\]
are open. 
After shrinking $Q$, we may assume that $\cL$ is torsion-free and $\cL^*$ is reflexive over any point of $Q$, $\cL$ is torsion-free, and $\cL^*$ is reflexive. 
Then $\cL^*\subset T_{\cU/Q}\subset T_\cU$ gives a family of foliations on $\cU$ over $Q$. 

Now for any $(X,\sF)\in\cS_P$, we dualize the following short exact sequence 
\[\xymatrix{0\ar[r] & \sF \ar[r] & T_X \ar[r]^-\tau & T_X/\sF \ar[r] & 0}\] 
and get the short exact sequence 
\[\xymatrix{0\ar[r] & (T_X/\sF)^* \ar[r]^-{\tau^*} & \Omega_X^{**} \ar[r] & \textnormal{coker}(\tau^*) \ar[r] & 0.}\] 
Note that 
\[\xymatrix{0\ar[r] & \textnormal{coker}(\tau^*) \ar[r] & \cO_X(K_\sF) \ar[r] & \cZ \ar[r] & 0}\] 
where $\cZ$ is the sheaf supported at a $0$-dimensional subscheme whose length is at most $\varphi$.  
Then we have 
\[\chi(X,\cO_X(K_\sF)(m\ell K_\sF+n\ell K_X)) - \chi(X,\textnormal{coker}(\tau^*)(m\ell K_\sF+n\ell K_X))\]
is a non-negative constant $\leq \varphi$. 
Thus, the function 
\[\chi(X,\textnormal{coker}(\tau^*)(m\ell K_\sF+n\ell K_X))\] 
has the form $\Phi(\ell)-s$. 
So by the construction of $Q$, there is a $q\in Q$ such that $\rho : X \rightarrow \cU_q$ is an isomorphism and 
$\textnormal{coker}(\tau^*)\cong\rho^*(\cL\vert_{\cU_q})$ and thus 
\[\sF \cong \textnormal{coker}(\tau^*)^*\cong\rho^*(\cL\vert_{\cU_q})^*.\] 

Notice that $\rho^*(\cL\vert_{\cU_q})^*$ and $\rho^*(\cL^*\vert_{\cU_q})$ are isomorphic at the generic point of $\cU_q$. 
By~\cite[Lemma 1.8]{hacon2021birational}, they are indeed isomorphic. 
Hence, $\cL^*$ gives the desired family of foliations. 
\qed
\end{enumerate}
\end{pf}

\section{Applications}
In this section, we show some applications of our boundedness result. 
First, we improve Proposition~\ref{bdd}, that is \cite[Proposition 4.1]{hacon2021birational}. 
\begin{prop}\label{dihedral_bdd}
Fix a function $P:\bZ_{\geq 0} \rightarrow \bZ$. 
There exist constants $C_1$ and $C_2$ such that, for any canonical model $(X_c,\sF_c)$ with Hilbert function $\chi(X_c,mK_{\sF_c})=P(m)$ for $m\geq 0$, the index (resp. embedding dimension) at each dihedral singularity 
is at most $C_1$ (resp. $C_2$). 
\end{prop}
\begin{pf}
By Theorem~\ref{key_thm}, there exists a family $\mu : \cX\rightarrow T$ of finite type 
and a family of foliations $\cF$ over $T$ such that 
for all foliated surfaces $(X,\sF)$ obtained as the minimal partial du Val resolution 
of some canonical foliation $(X_c,\sF_c)$ of general type with fixed Hilbert function $\chi(mK_{\sF_c}) = P(m)$, 
there exists a $t\in T$ and an isomorphism $\rho:X\rightarrow \cX_t$ such that $\sF\cong\rho^*(\cF\vert_{\cX_t})$. 
We may assume that 
\[T' = \{t\,\vert\, (\cX,\cF)_t \cong (X,\sF) \textnormal{ for some } (X,\sF)\in \cS_P\}\] 
is dense in $T$. 

After shrinking $T$, we may assume that $T$ is irreducible and $K_\cF$ is $\bQ$-Cartier. 
Let $D = i(K_\cF)K_\cF$ and fix a $t_0\in T'$. 
Since $D_{t_0}$ is nef, we have $D_{t}$ is nef for $t\in T^* = T\backslash\cup_{i=1}^\infty T_i$ where $T_i$'s are all proper closed subsets of $T$. 
(For a reference, see \cite[Proposition 1.4.14]{lazarsfeld2004positivity}.)
Note that, for all $t\in T^*$, we have $D_t^n = D_{t_0}^n>0$ since $D_{t_0}$ is nef and big.
Hence, $D_t$ is also nef and big for all $t\in T^*$. 

Let $S_N := \{t\in T\vert\, h^0(\cX_t,ND_{t})\geq 1\}$ which is a closed subset in $T$ 
by the upper semi-continuity for $N\in\bN$. 
Note that we have \[T(\bC) = \left(\bigcup_{N\in\bN} S_N(\bC)\right)\cup\left(\bigcup_i T_i(\bC)\right).\] 

By Baire category theorem, there is an $N$ such that $S_N = T$. 
Thus, the relative base locus of $D$ is proper. 
Let $U$ be an open subset of $T$ such that all irreducible components of the relative base locus of $D$ are flat over $U$. 
Notice that we may assume that $t_0\in U$. 

Let $\cE_j$'s be those irreducible components of the relative base locus of $D$ over $U$ of codimension one. 
Note that, for any $t\in T'\cap U$ and for any irreducible curve $E\subset \cY_t$ with $D\cdot E = 0$, we have $E$ is in the relative base locus of $D$ and hence $E\subset\cE_j$ for some $j$. 
Let $\pi:\wt{\cX}\rightarrow\cX$ be the log resolution of $\cX$ and the $\cE_j$'s with the exceptional divisors $\sum \cF_\ell$. 
By shrinking $T$ further, we may assume that $\nu := \mu\circ \pi$ is smooth and $\cF_\ell$ is horizontal with respect to $\nu$. 
Let $\wt{\cE_j}$ be the proper transform of $\cE_j$ via $\pi$. 

Note that $\Lambda := \{\wt{\cE_j}, \cF_\ell\}$ is a finite set. 
So there are only finitely many possible connected subsets of $\Lambda$ with the dual graph of $D_n$ type. 
By \cite[Satz 2.9 and Satz 2.11]{brieskorn1968rationale}, the intersection matrix of a dihedral singularity determines the embedding dimension and the group acting on it. 
Since the intersection matrix remains the same in the family, we have the boundedness of the embedding dimension and the order of the group acting on these singularities. 
Hence the index of $X$ at each dihedral singularities is also bounded.
\qed
\end{pf}

\begin{rmk}
\begin{enumerate}
\item It is possible that some connected subset of $\Lambda$ of $D_n$ type is redundant. 
\item By the same argument as above, we could also show the boundedness of the cyclic quotient singularities which are canonical non-terminal foliation singularities. 
Moreover, we could have the boundedness of the cusps. 
Precisely, over any cusp singularities, the intersection matrix and the number of the irreducible components of the exceptional divisor of the minimal resolution are bounded. 
\end{enumerate}
\end{rmk}

Also we answer \cite[Conjecture 1]{hacon2021birational}. 
\begin{thm}\label{bir}
Fix a function $P: \bZ_{\geq 0} \rightarrow \bZ$. 
There exists an integer $m_P$ such that if $(X,\sF)$ is a canonical model of a surface with $\kappa(K_\sF)=2$ and $\chi(mK_\sF)=P(m)$ for all $m\geq 0$, then for any $m>0$ divisible by $m_P$, $|mK_\sF|$ defines a birational map which is an isomorphism on the complement of the cusp singularities. 
\end{thm}
\begin{pf}
The proof follows closely from the proof of \cite[Theorem 4.3]{hacon2021birational}. 
We include the proof for the reader's convenience. 

Let $g : Y \rightarrow X$ be the minimal resolution of cusps and $\sG$ be the pullback foliation of $\sF$. 
Then $K_\sG = g^*K_\sF$ and $i(\sG)=i_\bQ(\sF)$. 
Since $K_\sF$ is numerically ample (see Lemma~\ref{ample}), we have that, for any curve $C$ on $X$,  $K_\sG\cdot g_*^{-1}C = K_\sF\cdot C>0$ and thus $i(\sG)K_\sG\cdot g_*^{-1}C\geq 1$. 
\begin{claim}
$K_Y+3i(\sG)K_\sG$ is nef. 
\end{claim}
\begin{pf}[Claim]
Given any irreducible curve $C$ on $Y$. We have following two cases. 
\begin{enumerate}
\item If $C$ is contracted by $g$, then 
\[(K_Y+3i(\sG)K_\sG)\cdot C = K_Y\cdot C = -2-C^2 \geq 0\] 
since $C^2\leq -2$ by the minimality of $g$. 
\item Suppose $C$ is not contracted by $g$. 
\begin{enumerate}
\item If $K_Y\cdot C\geq 0$, then $(K_Y+3i(\sG)K_\sG)\cdot C \geq 3i(\sG)K_{\sF}\cdot g_*C > 0$ by Lemma~\ref{ample}.

\item If $K_Y\cdot C<0$, by \cite[Proposition 3.8]{fujino2012minimal}, every $K_Y$-negative extremal ray is spanned by a rational curve $C$ with $0<-K_Y\cdot C\leq 3$. 
Thus, if we write $C \equiv \sum a_iC_i+D$ where $D\cdot K_Y\geq 0$, the $\bR_{\geq 0}[C_i]$'s are the $K_Y$-negative extremal rays, and the $a_i$'s are non-negative, then at least one of $a_i$ is positive. 
By what we have seen above, $(K_Y+3i(\sG)K_\sG)\cdot D\geq 0$. 
Thus, we have 
\begin{eqnarray*}
(K_Y+3i(\sG)K_\sG)\cdot C &=& \sum_ia_i (K_Y+3i(\sG)K_\sG)\cdot C_i + (K_Y+3i(\sG)K_\sG)\cdot D \\
&\geq & \sum_ia_i(- 3+3) + 0 = 0.
\end{eqnarray*}
\end{enumerate} 
\end{enumerate}
\end{pf}

Now, by \cite[Lemma 3.5]{hacon2021birational}, we know $\gamma K_\sG-K_Y$ is big for any fixed integer $\gamma$ with 
\[\gamma > \left\lceil\max\left\{\frac{2K_{\sF}\cdot K_{X}}{K_{\sF}^2}+3i_\bQ(\sF),0\right\}\right\rceil.\]

Note that $L := (\beta-\gamma)K_\sG+(\gamma K_\sG-K_Y)$ is pseudo-effective for any $\beta>\gamma$.
Let $L=P+N$ be the Zariski decomposition. 
Then we have 
\[P^2\geq(\beta-\gamma)^2K_\sG^2\geq \left(\frac{\beta-\gamma}{i_\bQ(\sF)}\right)^2.\]
And if $C$ is not contained in the exceptional divisor of $g$, then we have \[P\cdot C\geq (\beta-\gamma)K_\sG\cdot C\geq \frac{\beta-\gamma}{i_\bQ(\sF)}.\]

By Proposition~\ref{dihedral_bdd} and Proposition~\ref{bdd}, we have the boundedness of embedding dimension of each rational singularity on $Y$. 
Thus, we have a uniform bound $\delta$ such that $\delta_\zeta\leq\delta$ for $|\zeta|=2$. 
Now we fix a $\beta$ divisible by $i_\bQ(\sF)$ such that $P^2>\delta$ and $P\cdot C>\delta/2$ for any irreducible curve $C$. 
By \cite[Theorem 0.1]{langer2001adjoint} and Proposition~\ref{sep}, we have that $K_Y+L = \beta K_\sG$ is very ample on the complement of the exceptional divisor of $g$. 
By \cite[Theorem 6.2]{sakai1984weil}, $g_*\cO_Y(\beta K_\sG) = \cO_X(\beta K_\sF)$. 
Hence, $\beta K_\sF$ is very ample on the complement of the cusps. \qed
\end{pf}

\bibliographystyle{amsalpha}
\addcontentsline{toc}{chapter}{\bibname}
\bibliography{Fol_MPdVR}

\end{document}